\documentclass[11pt,a4paper,reqno]{amsart}  
\usepackage{amssymb} \usepackage{amscd} \usepackage{times}
\newcommand{\be} {\begin{equation}}
\newcommand{\ee} {\end{equation}}
\newcommand{\bea} {\begin{eqnarray}}
\newcommand{\eea} {\end{eqnarray}}
\newcommand{\Bea} {\begin{eqnarray*}}
\newcommand{\Eea} {\end{eqnarray*}}

\def\zbb{\mathbb{Z}}  
  
  \def\phi{\varphi}
 \def\p1{{\mathbb{P}^1_\zbb}}

\newtheorem{Theorem}{\quad Theorem}[section]

\begin{document}
\title{A compactness result for a Gelfand-Liouville system with Lipschitz condition.}
\author{Samy Skander Bahoura} 
\address{Departement de Mathematiques, Universite Pierre et Marie Curie, 2 place Jussieu, 75005, Paris, France.}
\email{samybahoura@yahoo.fr}
\maketitle
\begin{abstract}
We give  a quantization analysis to an elliptic system (Gelfand-Liouville type system) with Dirichlet condition. An application, we have a  compactness result for an elliptic system with  Lipschitz condition. 
\end{abstract}


\section{Introduction and Main Results}

We set $ \Delta = \partial_{11} + \partial_{22} $  on open set $ \Omega $ of $ {\mathbb R}^2 $ with a smooth boundary.

\bigskip

We consider the following equation:

$$ (P)   \left \{ \begin {split} 
       -\Delta u & = V e^{v} \,\, &\text{in} \,\, & \Omega  \subset {\mathbb R}^2, \\
      - \Delta v & = W e^{u} \,\, &\text{in} \,\, & \Omega  \subset {\mathbb R}^2, \\
                    u & = 0  \,\,             & \text{in} \,\,    &\partial \Omega,\\ 
                    v & = 0  \,\,             & \text{in} \,\,    &\partial \Omega.                        
   \end {split}\right.
$$

Here:

$$  0 \in \partial \Omega $$
 
When $ u=v $, the above system is reduced to an equation which was studied by many authors, with or without  the boundary condition, also for Riemann surfaces,  see [1-16],  one can find some existence and compactness results, also for a system.

Among other results, we  can see in [6] the following important Theorem,

\smallskip

{\bf Theorem A.}{\it (Brezis-Merle [6])}.{\it  Consider the case of one equation; if $ (u_i)_i=(v_i)_i $ and $ (V_i)_i=(W_i)_i $ are two sequences of functions relatively to the problem $ (P) $ with, $ 0 < a \leq V_i \leq b < + \infty $, then, for all compact set $ K $ of $ \Omega $,

$$ \sup_K u_i \leq c = c(a, b, K, \Omega). $$}

\smallskip

{\bf Theorem B} {\it (Brezis-Merle [6])}.{\it Consider the case of one equation and assume that $ (u_i)_i $ and $ (V_i)_i $ are two sequences of functions relatively to the previous problem $ (P) $ with, $ 0 \leq V_i \leq b < + \infty $, and,

$$ \int_{\Omega} e^{u_i} dy  \leq C, $$

then, for all compact set $ K $ of $ \Omega $,

$$ \sup_K u_i \leq c = c(b, C, K, \Omega). $$}

Next, we call energy the following quantity:

$$ E= \int_{\Omega} e^{u_i} dy. $$

The boundedness of the energy is a necessary condition to work on the problem $ (P) $ as showed in $ [6] $, by the following counterexample.

\smallskip

{\bf Theorem C} {\it (Brezis-Merle [6])}.{\it  Consider the case of one equation, then there are two sequences $ (u_i)_i $ and $ (V_i)_i $ of the problem $ (P) $ with, $ 0 \leq V_i \leq b < + \infty $, and,

$$ \int_{\Omega} e^{u_i} dy  \leq C, $$

and

$$ \sup_{\Omega}  u_i \to + \infty. $$}

Note that in [11], Dupaigne-Farina-Sirakov proved (by an existence result of Montenegro, see [14]) that the solutions of the above system when $ V $ and $ W $ are constants can be extremal and this condition imply the boundedness of the energy and directly the compactness.
Note that in [10], if we assume (in particular) that $ \nabla \log V $ and $ \nabla \log W $ and $ V >a >0 $ or $ W > a'>0 $ and $ V, W $ are nonegative and uniformly bounded then the energy is bounded and we have a compactness result.

Note that in the case of one equation, we can prove by using the Pohozaev identity that if $ + \infty >b \geq V \geq a >0 $, $ \nabla V $ is uniformely Lipschitzian that the energy is bounded when $ \Omega $ is starshaped. In [13] Ma-Wei, using the moving-plane method showed that this fact is true for all domain $ \Omega $ with the same assumptions on $ V $. In [10] De Figueiredo-do O-Ruf extend this fact to a system by using the moving-plane method for a system.

Theorem C, shows that we have not a global compactness to the previous problem with one equation, perhaps we need more information on $ V $ to conclude  to the boundedness of the solutions. When $ \nabla \log V $ is  Lipschitz function, Chen-Li and Ma-Wei see [7] and [13], showed that we have a compactness on all the open set. The proof is via the moving plane-Method of Serrin and Gidas-Ni-Nirenberg. Note that in [10], we have the same result for this system when $ \nabla \log V $ and $ \nabla \log W $ are uniformly bounded. We will see below that for a system we also have a compactness result when $ V $ and $ W $ are Lipschitzian.

Now consider the case of one equation. In this case our equation have nice properties.

If we assume $ V $ with more regularity, we can have another type of estimates, a $ \sup + \inf $ type inequalities. It was proved by Shafrir see [16], that, if $ (u_i)_i, (V_i)_i $ are two sequences of functions solutions of the previous equation without assumption on the boundary and, $ 0 < a \leq V_i \leq b < + \infty $, then we have the following interior estimate:

$$ C\left (\dfrac{a}{b} \right ) \sup_K u_i + \inf_{\Omega} u_i \leq c=c(a, b, K, \Omega). $$

\bigskip

Now, if we suppose $ (V_i)_i $ uniformly Lipschitzian with $ A $ the
Lipschitz constant, then, $ C(a/b)=1 $ and $ c=c(a, b, A, K, \Omega)
$, see [5]. 
\smallskip

Here we are interested by the case of a system of this type of equation. First, we give the behavior of the blow-up points on the boundary and in the second time we have a proof of  compactness of the solutions to Gelfand-Liouville type system with Lipschitz condition.

\smallskip

Here, we write an extention of Brezis-Merle Problem (see [6]) is:

\smallskip

{\bf Problem}. Suppose that $ V_i \to  V $ and  
$ W_i \to  W $ in $ C^0( \bar \Omega ) $, with, $ 0 \leq V_i \leq b_1  $ and $ 0 \leq W_i \leq b_2  $ for some positive constants $ b_1, b_2 $. Also, we consider a sequence of solutions $ (u_i), (v_i) $ of $ (P) $ relatively to $ (V_i), (W_i) $ such that,

$$ \int_{\Omega} e^{u_i} dx \leq C_1,\,\,\,  \int_{\Omega} e^{v_i} dx \leq C_2,  $$

is it possible to have:

$$ ||u_i||_{L^{\infty}}\leq C_3=C_3(b_1, b_2, C_1, C_2, \Omega) ? $$

and,

$$ ||v_i||_{L^{\infty}}\leq C_4=C_4(b_1, b_2, C_1, C_2, \Omega) ? $$

In this paper we give a caracterization of the behavior of the blow-up points on the boundary and also a proof of the compactness theorem when $ V_i $ and $ W_i $ are uniformly Lipschitzian. For the behavior of the blow-up points on the boundary, the following condition are enough,

$$ 0 \leq  V_i \leq b_1, \,\,\, 0 \leq  W_i \leq b_2, $$

The conditions $ V_i \to  V $ and $ W_i \to  W $ in $ C^0(\bar \Omega) $ are not necessary.

\bigskip

But for the proof of the compactness  for the  Gelfand-Liouville type system (Brezis-Merle type problem) we assume that:

$$ ||\nabla V_i||_{L^{\infty}}\leq  A_1, \,\,\, ||\nabla W_i||_{L^{\infty}}\leq  A_2. $$

Our main result are:

\begin{Theorem}  Assume that $ \max_{\Omega} u_i \to +\infty $ and $ \max_{\Omega} v_i \to +\infty $ Where $ (u_i) $ and $ (v_i) $ are solutions of the probleme $ (P) $ with:
 
 $$ 0 \leq V_i \leq b_1,\,\,\, {\rm and } \,\,\, \int_{\Omega}  e^{u_i} dx \leq C_1, \,\,\, \forall \,\, i, $$
 
and,

 $$ 0 \leq W_i \leq b_2,\,\,\, {\rm and } \,\,\, \int_{\Omega}  e^{v_i} dx \leq C_2, \,\,\, \forall \,\, i, $$

 then;  after passing to a subsequence, there is a finction $ u $,  there is a number $ N \in {\mathbb N} $ and $ N  $ points $ x_1, x_2, \ldots, x_N \in  \partial \Omega $, such that, 

$$ \int_{\partial \Omega} \partial_{\nu} u_i  \phi  \to \int_{\partial \Omega} \partial_{\nu} u  \phi +\sum_{j=1}^N \alpha_j  \phi(x_j), \, \alpha_j \geq 4\pi, $$
$ \text{for\, any}\,\, \phi\in C^0(\partial \Omega) $, 
and,

$$ u_i \to u \,\,\, {\rm in }\,\,\, C^1_{loc}(\bar \Omega-\{x_1,\ldots, x_N \}). $$

$$ \int_{\partial \Omega} \partial_{\nu} u_i  \phi  \to \int_{\partial \Omega} \partial_{\nu} u  \phi +\sum_{j=1}^N \beta_j  \phi(x_j), \, \beta_j \geq 4\pi, $$
$ \text{for\, any}\,\, \phi\in C^0(\partial \Omega) $,  and,

$$ v_i \to v \,\,\, {\rm in }\,\,\, C^1_{loc}(\bar \Omega-\{x_1,\ldots, x_N \}). $$

\end{Theorem} 

 In the following theorem, we have a proof for the global a priori estimate which concern the problem $ (P) $.

\bigskip

\begin{Theorem}Assume that $ (u_i), (v_i) $ are solutions of $ (P) $ relatively to $ (V_i), (W_i) $ with the following conditions:

$$ x_1=0 \in \partial \Omega, $$

and,

$$ 0 \leq  V_i \leq b_1, \,\,  ||\nabla V_i||_{L^{\infty}} \leq A_1,\,\, {\rm and } \,\,\, \int_{\Omega} e^{u_i} \leq C_1, $$

$$ 0 \leq  W_i \leq b_2, \,\,  ||\nabla W_i||_{L^{\infty}} \leq A_2,\,\, {\rm and } \,\,\, \int_{\Omega} e^{v_i} \leq C_2, $$

We have,

$$  || u_i||_{L^{\infty}} \leq C_3(b_1, b_2, A_1, A_2, C_1, C_2, \Omega), $$

and,

$$  || v_i||_{L^{\infty}} \leq C_4(b_1, b_2, A_1, A_2, C_1, C_2, \Omega), $$

\end{Theorem} 

\section{Proof of the theorems} 

\bigskip

\underbar {\it Proof of theorem 1.1:} 

\bigskip

Since $ V_ie^{v_i} $ and $ W_ie^{u_i} $ are bounded in $ L^1(\Omega) $, we can extract from those two sequences two subsequences which converge to two nonegative measures $ \mu_1 $ and $ \mu_2 $.

If $ \mu_1(x_0) < 4 \pi $, by a Brezis-Merle estimate for the first equation, we have $ e^{u_i} \in L^{1+\epsilon} $ around $ x_0 $, by the elliptic estimates, for the second equation, we have $ v_i \in W^{2, 1+\epsilon} \subset L^{\infty} $ around $ x_0 $, and , returning to the first equation, we have $ u_i \in L^{\infty} $ around $ x_0 $.

If $ \mu_2(x_0) < 4 \pi $, then $ u_i $ and $ v_i $ are also locally bounded around $ x_0 $.

Thus, we take a look to the case when, $ \mu_1(x_0) \geq 4 \pi $ and $ \mu_2(x_0) \geq 4 \pi $. By our hypothesis, those points $ x_0 $ are finite.

We will see that inside $ \Omega $ no such points exist. By contradiction, assume that, we have $ \mu_1(x_0) \geq 4 \pi $. Let us consider a ball $ B_R(x_0) $ which contain only $ x_0 $ as nonregular point. Thus, on $ \partial B_R(x_0) $, the two sequence $ u_i $ and $ v_i $ are uniformly bounded. Let us consider:

$$ \left \{ \begin {split} 
       -\Delta z_i & = V_i e^{v_i} \,\, &\text{in} \,\, & B_R(x_0)  \subset {\mathbb R}^2, \\
                       z_i & = 0  \,\,             & \text{in} \,\,    &\partial B_R(x_0).                                            
   \end {split}\right.
$$

By the maximum principle we have:

$$ z_i \leq u_i $$ 

and $ z_i \to z $ almost everywhere on this ball, and thus,

$$ \int e^{z_i} \leq \int e^{u_i} \leq C, $$

and,

$$ \int e^z \leq C.$$

but, $ z  $ is  a solution to the following equation:

$$ \left \{ \begin {split} 
       -\Delta z & = \mu_1\,\, &\text{in} \,\, & B_R(x_0) \subset {\mathbb R}^2, \\
                     z & = 0  \,\,             & \text{in} \,\,    &\partial B_R(x_0).                                            
   \end {split}\right.
$$

with, $ \mu_1 \geq 4 \pi $ and thus, $ \mu_1 \geq 4\pi \delta_{x_0} $ and then, by the maximum principle:

$$ z \geq -2 \log |x-x_0|+ C $$

thus,

$$ \int e^z = + \infty, $$

which is a contradiction. Thus, there is no nonregular points inside  $ \Omega $

Thus, we consider the case where we have nonregular points on the boundary, we use two estimates:

$$ \int_{\partial \Omega} \partial_{\nu} u_i d\sigma \leq C_1,\,\,\, \int_{\partial \Omega} \partial_{\nu} v_i d\sigma \leq C_2, $$
 
 and,
 
 $$ ||\nabla u_i||_{L^q} \leq C_q, \,\,\, ||\nabla v_i||_{L^q} \leq C'_q, \,\,\forall \,\, i\,\, {\rm and  }  \,\, 1< q < 2. $$

We have the same computations, as in the case of one equation.

We consider a points $ x_0 \in \partial \Omega $ such that:

$$ \mu_1(x_0) < 4 \pi. $$

We consider a test function on the boundary $ \eta $ we extend $ \eta $ by a harmonic function on $ \Omega $, we write the equation:

$$ -\Delta ((u_i-u)\eta) =(V_i e^{v_i}-Ve^v)\eta+ <\nabla (u_i-u)|\nabla \eta> = f_i $$

with,

$$ \int |f_i| \leq 4 \pi-\epsilon +o(1) < 4\pi-2\epsilon <4\pi, $$

$$ -\Delta ((v_i-v)\eta) =(W_i e^{u_i}-We^u)\eta+ <\nabla (v_i-v)|\nabla \eta> = g_i, $$

with,

$$ \int |g_i| \leq 4 \pi-\epsilon +o(1) < 4\pi-2\epsilon <4\pi, $$

By the Brezis-Merle estimate, we have uniformly, $ e^{u_i} \in L^{1+\epsilon} $ around $ x_0 $, by the elliptic estimates, for the second equation,  we have $ v_i \in W^{2, 1+\epsilon} \subset L^{\infty} $ around $ x_0 $, and , returning to the first equation, we have $ u_i \in L^{\infty} $ around $ x_0 $.

We have the same thing if we assume:

$$ \mu_2(x_0) < 4 \pi. $$

Thus, if $ \mu_1(x_0) < 4 \pi $ or $ \mu_2(x_0) < 4 \pi $, we have for $ R >0 $ small enough:

$$ (u_i,v_i) \in L^{\infty}(B_R(x_0)\cap \bar \Omega). $$

By our hypothesis the set of the points such that:

$$ \mu_1(x_0)  \geq  4 \pi,  \,\,\, \mu_2(x_0)   \geq 4 \pi, $$

is finite, and, outside this set $ u_i $ and $ v_i $ are locally uniformly bounded. By the elliptic estimates, we have the $ C^1 $ convergence to  $ u $ and $ v $  on each compact set of $ \bar \Omega- \{x_1, \ldots x_N\} $.

\bigskip

\underbar {\it Proof of theorem 1.2:} 

\bigskip

Without loss of generality, we can assume that $ 0 $ is a blow-up point (either, we use a translation). Also, by a conformal transformation, we can assume that $ \Omega =B_1^+ $, the half ball, and $ \partial^+ B_1^+ $ is the exterior part, a part which not contain $ 0 $ and on which  $ u_i $ and $ v_i $ converge in the $ C^1 $ norm to $ u $ and $ v $. Let us consider $ B_{\epsilon}^+ $, the half ball with radius $ \epsilon >0 $.

\bigskip 

The Pohozaev identity gives :

\be  \int_{B_{\epsilon}^+} \Delta u_i < x |\nabla v_i > dx = -  \int_{B_{\epsilon}^+} \Delta v_i < x |\nabla u_i > dx + \int_{\partial^+ B_{\epsilon}^+}  g(\partial_{\nu}u_i, \partial_{\nu}v_i)d\sigma, \label{(10)}\ee

Thus,

\be  \int_{B_{\epsilon}^+} V_i e^{v_i}< x |\nabla v_i > dx = -  \int_{B_{\epsilon}^+} W_ie^{u_i} < x |\nabla u_i > dx + \int_{\partial^+ B_{\epsilon}^+}  g(\partial_{\nu}u_i, \partial_{\nu}v_i)d\sigma, \label{(10)}\ee

After integration by parts, we obtain:

$$  \int_{B_{\epsilon}^+} V_i e^{v_i} dx +  \int_{B_{\epsilon}^+} < x |\nabla V_i > e^{v_i} dx+  \int_{\partial B_{\epsilon}^+} < \nu |\nabla V_i > d\sigma+ $$

$$ + \int_{B_{\epsilon}^+} W_i e^{u_i} dx +  \int_{B_{\epsilon}^+} < x |\nabla W_i > e^{u_i} dx+  \int_{\partial B_{\epsilon}^+} < \nu |\nabla W_i > d\sigma = $$

$$ =  \int_{\partial^+ B_{\epsilon}^+} g(\partial_{\nu}u_i, \partial_{\nu}v_i)d\sigma, \label{(10)} $$ 

Also, for $ u $ and $ v $, we have:

$$  \int_{B_{\epsilon}^+} Ve^{v} dx +  \int_{B_{\epsilon}^+} < x |\nabla V > e^{v} dx+  \int_{\partial B_{\epsilon}^+} < \nu |\nabla V> d\sigma+ $$

$$ + \int_{B_{\epsilon}^+} We^{u} dx +  \int_{B_{\epsilon}^+} < x |\nabla W> e^{u} dx+  \int_{\partial B_{\epsilon}^+} < \nu |\nabla W > d\sigma = $$

$$ =  \int_{\partial^+ B_{\epsilon}^+} g(\partial_{\nu}u, \partial_{\nu}v)d\sigma, $$ 

If, we take the difference, we obtain:

$$ (1+o(\epsilon))(\int_{B_{\epsilon}^+} V_i e^{v_i} dx -\int_{B_{\epsilon}^+} V e^{v} dx) + $$

$$ + (1+o(\epsilon))(\int_{B_{\epsilon}^+} W_i e^{u_i} dx -\int_{B_{\epsilon}^+} W e^{u} dx) = $$

$$ = \alpha_1+\beta_1+ o(\epsilon)+o(1)= o(1), $$

a contradiction.

\end{document}